\documentclass[11pt,leqno,fleqn]{article}
\usepackage{epsf,amsmath,amssymb,amscd,hyperref,srcltx}
\usepackage{amsfonts}
\usepackage{counters}
\usepackage{psfrag}
\newcommand{\dps}{\displaystyle}

\newcommand{\kies}[2]{\mbox{${{#1}\choose{#2}}$}}

\newcounter{figuren}

\newcommand{\dy}[2]{%
\refstepcounter{equation}%
\LABEL{#1}%
\begin{list}{}{
\topsep 5mm
\leftmargin 18mm
\rightmargin 0cm
\itemsep 0mm
\listparindent 0mm
\parsep 0mm
\itemsep 0mm
\labelsep 0mm
\labelwidth 18mm
}%
\item[\rm (\theequation)\hfill]
#2
\end{list}%
}
\newcommand{\dyz}[1]{%
\refstepcounter{equation}%
\begin{list}{}{
\topsep 5mm
\leftmargin 18mm
\rightmargin 0cm
\itemsep 0mm
\listparindent 0mm
\parsep 0mm
\itemsep 0mm
\labelsep 0mm
\labelwidth 18mm
}%
\item[\rm (\theequation)\hfill]
#1
\end{list}%
}
\newcommand{\dyyz}[1]{\dyz{\raggedright$\dps#1$}}

\newcommand{\de}[2]{\dy{#1}{\raggedright$\displaystyle #2 $}}
\newcommand{\dez}[1]{\dyz{\raggedright$\displaystyle #1 $}}
\newcommand{\leeg}[1]{}

\newcounter{stelling}
\newcommand{\thm}[2]{\setcounter{gevolg}{0}\setcounter{claim}{0}\refstepcounter{stelling}\vspace{4mm}\noindent{\bf Theorem \thestelling.}\label{#1}{\it #2}}

\newalphacounterwithin{gevolg}{stelling}
\newcommand{\cor}[2]{\refstepcounter{gevolg}\setcounter{claim}{0}\vspace{4mm}\noindent{\bf Corollary \thegevolg.}\label{#1}{\it #2}}

\newcounter{hulpstelling}

\newcounter{bewering}

\newcounter{claim}

\newalphacounterwithin{subclaim}{claim}

\newcounter{opmerking}

\newcounter{hoofdstuk}

\newcounter{sectie}
\newcounter{subsectie}
\newcounterwithin{ex}{sectie}

\newcommand{\sectz}[1]{\refstepcounter{sectie}\setcounter{subsectie}{0}\setcounter{ex}{0}
\section*{\boldmath \thesectie. #1}%
}

\newcounter{lit}

\newcommand{\pf}{\vspace{3mm}\noindent{\bf Proof.}\ }

\newcommand{\bx}{\hspace*{\fill} \hbox{\hskip 1pt \vrule width 4pt height 8pt depth 1.5pt \hskip 1pt}

\addvspace{4mm}}

\newcommand{\rf}[1]{{\rm (\ref{#1})}}

\newcommand{\kint}[2]{\mbox{$\int$}}

\newcommand{\NIET}[1]{}

\newcommand{\LABEL}[1]{\label{#1}}

\newcommand{\rank}{\text{\rm rank}}

\newcommand{\oN}{{\mathbb{N}}}

\newcommand{\oR}{{\mathbb{R}}}

\newcommand{\op}{_{\text{\rm o}}}
\begin{document}

\begin{center}
{\LARGE\bf Low rank approximation of polynomials

}
\vspace{4mm}

{\large
\hspace{10mm}
Alexander Schrijver\footnote{ CWI and University of Amsterdam.
Mailing address: CWI, Kruislaan 413, 1098 SJ Amsterdam,
The Netherlands.
Email: lex@cwi.nl.}}

\end{center}

\noindent
{\small{\bf Abstract.}
Let $k\leq n$.
Each polynomial $p\in\oR[x_1,\ldots,x_n]$ can be uniquely written as
$p=\sum_{\mu}\mu p_{\mu}$, where $\mu$ ranges over the set $M$ of all monomials in $\oR[x_1,\ldots,x_k]$
and where $p_{\mu}\in\oR[x_{k+1},\ldots,x_n]$.
If $p$ is $d$-homogeneous and $\varepsilon>0$, we say
that $p$ is {\em $\varepsilon$-concentrated on the first $k$ variables} if
$$
\sum_{\mu\in M\atop\deg(\mu)<d}\max_{x\in\oR^{n-k}\atop\|x\|=1}p_{\mu}(x)^2
\leq
\varepsilon\|p\|^2,
$$
where $\|p\|$ is the Bombieri norm of $p$.
We show
that for each $d\in\oN$ and $\varepsilon>0$ there exists
$k_{d,\varepsilon}$ such that for each $n$ and each $d$-homogeneous
$p\in\oR[x_1,\ldots,x_n]$
there exists $k\leq k_{d,\varepsilon}$
such that $p$ is $\varepsilon$-concentrated on the first $k$ variables
{\em after some orthogonal transformation of $\oR^n$}.
(So $k_{d,\varepsilon}$ is independent of the number $n$ of variables.)

We derive this as a consequence of a more general theorem on low rank
approximation of polynomials.

}

\sectz{Introduction}

A polynomial is said to have {\em rank} 1 if it is nonzero and a product of linear functions.
The {\em rank} of a polynomial $p$, denoted by $\rank(p)$, is the minimum number of
rank 1 polynomials that add up to $p$.
Low rank of polynomials helps in computing them.
Therefore, we investigate in how much polynomials can be approximated by low rank polynomials.
In particular, we consider $d$-homogeneous polynomials
(i.e., homogeneous polynomials of total degree $d$).

In fact there are a number of different notions of rank function of a polynomial, but for what
follows we can take any of them, as long as the function is invariant under orthogonal
transformations of the polynomial and it does not increase if we set variables in the polynomial to 0.
So one can also take the minimum number of linear functions that generate $p$, or
the {\em Waring rank}: the minimum number of powers of linear functions that linearly span $p$.
Also one may take the minimum number of rank 1 polynomials of Bombieri
norm at most 1 that add up to $p$ (for the definition of Bombieri norm see below).
More strongly, one can take any monotonically nondecreasing function of
any rank function.

We recall the
{\em Bombieri norm} $\|p\|$ of a $d$-homogeneous polynomials
$p\in\oR[x_1,\ldots,x_n]$:
\dez{
\|p\|:=
(\sum_{\alpha\in\oN^n}\kies{d}{\alpha_1,\ldots,\alpha_n}^{-1}p_{\alpha}^2)^{1/2},
}
where, for $\alpha\in\oN^n$, $p_{\alpha}$ be the coefficient of
$x_1^{\alpha_1}\cdots x_n^{\alpha_n}$.\footnote{The Bombieri norm is the tensor norm if we view $d$-homogeneous polynomials
as symmetric $d$-tensors.
In our estimates we can replace the Bombieri norm equivalently by
the square root of the sum of the squares of the coefficients, as
we fix $d$ and as
the two values are bounded by each other up to a factor of $d!$.
However, the Bombieri norm behaves better algebraically --- in particular,
it is invariant under orthogonal transformations of $\oR^n$.}

Let us remark that it is not true that for each $d\in\oN$ and $\varepsilon>0$
there is a $k\in\oN$ such that
for each $d$-homogeneous polynomial $p$ there exists a $d$-homogeneous polynomial $q$ of rank
at most $k$ such that $\|p-q\|\leq\varepsilon\|p\|$.
This is shown by the polynomials $\sum_{i=1}^nx_i^2$.

However, the approximation becomes valid by considering the {\em operator norm}:
\dez{
\|p\|\op:=\max_{x\in\oR^n\atop\|x\|=1}|p(x)|.
}
Then Fernandez de la Vega, Kannan, Karpinski, and Vempala [1] showed
(it is also a special case of the `weak regularity for Hilbert spaces'
in Lov\'asz and Szegedy [2])
\dy{10ok12c}{
for each $d\in\oN$ and $\varepsilon>0$ there exists $k\in\oN$ such that for each
$d$-homogeneous polynomial $p$ (in any number of variables) there exists a
$d$-homogeneous polynomial $q$ of rank $\leq k$ such that
$\|p-q\|\op\leq\varepsilon\|p\|$ for each $x$.
}
Important in \rf{10ok12c} is that $k$ is independent of the number of variables.
In fact, in [1] it is shown
that one can take $k=\lfloor\varepsilon^{-2}\rfloor$.

In this paper we give an extension of \rf{10ok12c}, using the compactness
result of [3] as main tool.

\sectz{Stronger low rank approximation}

To describe the extension,
let for any $n$-variable polynomial $p$ and any subspace $V$ of $\oR^n$,
$\pi_V$ denote the orthogonal projection onto $V$, and $p_V:=p\circ\pi_V$.
Then $\|p\|\op$ may be alternatively described as
the supremum of $\|p_V\|$ taken over all 1-dimensional subspaces of $\oR^n$.
As an extension, define for any $k\in\oN$:
\dez{
\|p\|_{(k)}:=\sup_{V\text{ subspace}\atop\dim V\leq k}\|p_V\|.
}
So $\|p\|\op=\|p\|_{(1)}$.
Moreover, as $\|p_V\|\leq\|p\|$ for any subspace $V$, we have $\|p\|_{(k)}\leq\|p\|$.

Let $P$ be the set of all $d$-homogeneous polynomials in $\oR[x_1,x_2,\ldots]$, each using
only a finite number of variables.
Call a function $r:P\to\oR$ {\em monotone} if it is invariant under
orthogonal transformations of the space and if it does not increase if
we set variables in a polynomial to 0.
Each of the above rank functions is monotone, and also each of the above
norms.
Moreover, if $r$ is monotone and $f:\oR\to\oR$ is monotonically increasing,
then also $f\circ r$ is monotone.

Then in \rf{10ok12c}, one may replace $\|p-q\|\op$
by $\|p-q\|_{f(\rank(q))}$, for any fixed function $f:\oN\to\oN$
given in advance.
In fact:

\thm{3ok12b}{
For each $d\in\oN$, $\varepsilon>0$, and monotone $r:P\to\oN$,
there exists $k_{d,\varepsilon,r}$ such that for each $p\in P$
there exists $q\in P$ with $r(q)\leq k_{d,\varepsilon,r}$
and $\|p-q\|_{(r(q))}\leq\varepsilon\|p\|$.
}

\pf
Let $d\in\oN$, $\varepsilon>0$, and $r:P\to\oN$ be given.
As $P$ has an inner product associated with the Bombieri norm $\|.\|$,
the completion $H$ of $P$ is a Hilbert space.
As for each $k$, $\|.\|_{(k)}\leq\|.\|$, the norm $\|.\|_{(k)}$ extends uniquely to $H$.
Let $B(P)$ and $B(H)$ be the closed unit balls of $P$ and $H$ respectively.

Then for each fixed $k$,
the norm $\|.\|_{(k)}$ is continuous with respect to the $\|.\|\op$-topology on $H$.
This follows from the fact that there exists a $c$ such that
$\|p\|_{(k)}\leq c\|p\|\op$ for each $p\in P$.
Indeed, consider any subspace $V$ of dimension $k$, which we may assume to be
$\oR^k$.
Now the collection $C$ of all $d$-homogeneous polynomials $p\in\oR^k$ with
$\|p\|\op\leq 1$ is bounded.
Otherwise, as $C$ is a compact convex set in finite dimensions, there would
be a nonzero $p$ with $\|p\|\op=0$.
However, if $p\neq 0$, then $p(x)\neq 0$ for some $x$.
Hence there exists a $c$ such that $\|p\|\leq c$ for all $p\in C$, which is $c$ as required.

Define, for each $q\in P$,
\dez{
U_q:=\{p\in B(H)\mid \|p-q\|_{(r(q))}<\varepsilon\}.
}
So $U_q$ is open in the $\|.\|\op$-topology.
Moreover, the $U_q$ for $q\in B(P)$ cover $B(H)$.
Indeed, for any $p\in B(H)$ there exists $q\in B(P)$ with
$\|p-q\|<\varepsilon$.
Then $\|p-q\|_{(r(q))}\leq\|p-q\|<\varepsilon$, so $p\in U_q$.

Let $G$ be the group of all transformations that consist of an orthogonal
transformation of $\oR^n$ for some $n$, leaving the other coordinates
invariant.
Then $G$ acts naturally on $H$.
By [3], the orbit space $(B(H),\|.\|\op)/G$ is compact.
(This is the quotient topological space of the topological space
$(B(H),\|.\|\op)$, taking the orbits of $G$ as quotient classes.)
Hence there is a finite set $Q\subseteq B(P)$
such that voor each $p\in B(P)$ there exist $q\in Q$ and $\psi\in G$ such that
$p^{\psi}\in U_q$.
Let $k_{d,\varepsilon,f}:=\max\{r(q)\mid q\in Q\}$.
We show that $k_{d,\varepsilon,f}$ is as required.

Let $p\in P$, say $p\in\oR[x_1,\ldots,x_n]$.
We may assume that $\|p\|=1$.
Then there exist $q\in Q$ and $\psi\in G$ such that $p^{\psi}\in U_q$.
So $\|p^{\psi}-q\|_{(r(q))}<\varepsilon$.
As $\|.\|_{(r(q))}$ is $G$-invariant, this gives, setting $\varphi:=\psi^{-1}$,
$\|p-q^{\varphi}\|_{(r(q))}<\varepsilon$.
Let $q'$ be the orthogonal projection of $q^{\varphi}$ onto $\oR[x_1,\ldots,x_n]$.
So $r(q')\leq r(q^{\varphi})=r(q)\leq k_{d,\varepsilon,f}$.
Then $\|p-q'\|_{(r(q'))}<\varepsilon$, for
let $U$ be a subspace of $\oR^n$ of dimension $\leq r(q')$.
Then $(q')_U=(q^{\varphi})_U$, and hence $\|(p-q')_U\|=\|(p-q^{\varphi})_U\|\leq
\|p-q^{\varphi}\|_{(r(q'))}<\varepsilon$.
\bx

\sectz{$\varepsilon$-concentration on first $k$ variables}

Let $p\in\oR[x_1,\ldots,x_n]$ be $d$-homogeneous, and let $k\leq n$.
For $\alpha\in\oN^k$, denote $x^{\alpha}:=x_1^{\alpha_1}\cdots x_k^{\alpha_k}$.
Then $p$ can be uniquely written as
\dez{
p=\sum_{\alpha\in\oN^k\atop |\alpha|=d}x^{\alpha}p_{\alpha},
}
where $p_{\alpha}\in\oR^{n-k}$ and $|\alpha|:=\alpha_1+\cdots+\alpha_k$.
We say that $p$ is {\em $\varepsilon$-concentrated on $x_1,\ldots,x_k$} if
\de{14ok12a}{
\sum_{\alpha\in\oN^k\atop|\alpha|<d}
\|p_{\alpha}\|\op^2
\leq
\varepsilon\|p\|.
}

\cor{7ok12e}{
For each $d$ and $\varepsilon>0$ there exists $k_{d,\varepsilon}$
such that for each $n$ and each $d$-homogeneous $p\in\oR[x_1,\ldots,x_n]$
there exist $k\leq k_{d,\varepsilon}$ and an orthogonal transformation
$\varphi$ of $\oR^n$ such that $p^{\varphi}$ is $\varepsilon$-concentrated
on the first $k$ variables.
}

\pf
Let $f:\oN\to\oN$ be defined by $f(k):=k+$ the number of $k$-variable
monomials of degree $<d$.
For any $d$-homogeneous polynomial $q$, let $w(q)$ be its Waring rank
(the minimum number of $d$-powers of linear functions that span $q$), and
Let $r(q):=f(w(q))$.
Let $k_{d,\varepsilon}:=k_{d,\varepsilon/d!,r}$, where the latter is taken from
Theorem \ref{3ok12b}.
We show that $k_{d,\varepsilon}$ is as required.

Let $p\in\oR[x_1,\ldots,x_n]$, and let $q\in P$ be as given by Theorem \ref{3ok12b}.
Then $k:=w(q)\leq r(q)\leq k_{d,\varepsilon}=k_{d,\varepsilon/d!,r}$.
As $q$ has Waring rank $k$,
there exists a $k$-dimensional subspace $U$ of $\oR^n$ such that $q=q_U$.
By applying an orthogonal transformation of $\oR^n$,
we can assume that $U=\oR^k$.
So $q\in\oR[x_1,\ldots,x_k]$.
We prove that $p$ satisfies \rf{14ok12a}.

For each $\alpha\in\oN^k$ with $|\alpha|<d$,
choose $z_{\alpha}\in\oR^{n-k}$ with $\|z_{\alpha}\|=1$ maximizing
$|p_{\alpha}(z_{\alpha})|$.
Let $V$ be the space spanned by $\oR^k$ and by the $z_{\alpha}$.
Then $\|p_{\alpha}\|\op=
\|p_{\alpha}(z_{\alpha})|\leq|(p_V)_{\alpha}\|=\|(p_V-p_U)_{\alpha}\|$ for each such $\alpha$.
Moreover, $\dim(V)\leq f(k)=f(w(q))=r(q)$.
So
\dyyz{
\sum_{\alpha\in\oN^k\atop|\alpha|<d}
\|p_{\alpha}\|\op^2
\leq
\sum_{\alpha\in\oN^k\atop|\alpha|\leq d}
\|(p_V-p_U)_{\alpha}\|^2
=
\sum_{\alpha\in\oN^k\atop|\alpha|\leq d}
\sum_{\beta\in\oN^{n-k}\atop|\beta|=d-|\alpha|}
\kies{d-|\alpha|}{\beta_1,\ldots,\beta_{n-k}}^{-1}
(p_V-p_U)_{\alpha,\beta}^2
\leq
d!
\sum_{\alpha\in\oN^k\atop|\alpha|\leq d}
\kies{d}{\alpha_1,\ldots,\alpha_k}^{-1}
\sum_{\beta\in\oN^{n-k}\atop|\beta|=d-|\alpha|}
\kies{d-|\alpha|}{\beta_1,\ldots,\beta_{n-k}}^{-1}
(p_V-p_U)_{\alpha,\beta}^2
=
d!
\|p_V-p_U\|^2
\leq
d!
\|p_V-q\|^2
=
d!
\|(p-q)_V\|^2
\leq
d!
\|p-q\|_{(r(q))}^2
\leq
\varepsilon^2\|p\|^2.
}
Here $\|p_V-p_U\|\leq \|p_V-q\|$ follows from the fact that
$q=q_U$ and $p_U=(p_V)_U$, so $p_U$ is the polynomial defined
on $U$ closest to $p_V$.
\bx

Let $\|p\|_{\infty}$ be the maximum absolute value of the coefficients of
$p$.
(This norm is not invariant under orthogonal transformations.)
We note that for any $d$ there is a $c$ such that $\|p\|_{\infty}\leq c\|p\|$
for each $d$-homogeneous polynomial.
This follows from the facts that $\|p\|_{\infty}=\|p_V\|_{\infty}$ for some
$d$-dimensional subspace $V$ of $\oR^n$ and that the set $C$ of $d$-homogeneous
polynomials in $\oR[x_1,\ldots,x_d]$ with $\|p\|\leq 1$ is bounded.
This implies that in Corollary \ref{7ok12e} one may replace $\|p_{\alpha}\|\op$
by $\|p_{\alpha}\|_{\infty}$.

\section*{References}\label{REF}
{\small
\begin{itemize}{}{
\setlength{\labelwidth}{4mm}
\setlength{\parsep}{0mm}
\setlength{\itemsep}{1mm}
\setlength{\leftmargin}{5mm}
\setlength{\labelsep}{1mm}
}
\item[\mbox{\rm[1]}] W. Fernandez de la Vega, R. Kannan, M. Karpinski, S. Vempala, 
Tensor decomposition and approximation schemes for constraint satisfaction problems,
in: {\em Proceedings of the 37th Annual {ACM} Symposium on Theory of Computing}
({STOC}'05),
pp. 747--754,
{ACM}, New York, 2005.

\item[\mbox{\rm[2]}] L. Lov\'asz, B. Szegedy, 
Szemer\'edi's lemma for the analyst,
{\em Geometric and Functional Analysis} 17 (2007) 252--270.

\item[\mbox{\rm[3]}] G. Regts, A. Schrijver, 
Compact orbit spaces in Hilbert spaces and limits of edge-colouring models,
preprint, 2012.
ArXiv \url{http://arxiv.org/abs/1210.2204}

\end{itemize}
}

\end{document}